\documentclass[a4paper,12pt]{article}
\usepackage[english]{babel}
\usepackage[utf8x]{inputenc}
\usepackage{layout}
\usepackage{color}
\usepackage[dvipsnames]{xcolor}
\usepackage{lipsum}
\usepackage{amsmath}
\usepackage{graphicx}
\usepackage{eurosym}
\usepackage[font=footnotesize, margin=1cm]{caption}
\usepackage{subcaption}
\usepackage{graphicx}
\usepackage{amssymb}
\usepackage{tikz}
\usepackage{fullpage}
\usepackage{nicefrac}

\usepackage{lineno}

\usepackage{hyperref}
\hypersetup{colorlinks,bookmarksopen,bookmarksnumbered,
linkcolor=MidnightBlue,pdfstartview=FitH,urlcolor=black,citecolor=MidnightBlue}

\title{{\color{MidnightBlue}{\bfseries Can $\pi$ generate itself? A Monte Carlo analysis of 314 trillion digits}}}
\author{
Alessandro Razeto\textsuperscript{a} and Nicola Rossi\textsuperscript{a} \\
  {\small \textsuperscript{a} Laboratori Nazionali del Gran Sasso (INFN), Via G. Acitelli, 22 67100 Assergi L'Aquila, Italy } 
}

\date{\today}

\begin{document}
\maketitle

\begin{abstract}
At the end of 2025, a record computation of $\pi$ reached 314 trillion decimal digits, providing the largest numerical dataset ever generated for this constant. We exploit this unprecedented dataset to investigate whether the digits of $\pi$ themselves can serve as a source of pseudorandom numbers for estimating $\pi$ through the simplest Monte Carlo method. Our results go beyond the normality hypothesis by providing empirical evidence of a high degree of statistical randomness in the available digits, although not of digit independence, which cannot hold for a deterministic sequence. By optimizing the mapping of the digit sequence into Monte Carlo samples, we obtain the highest precision allowed by the dataset. As predicted, the method successfully reproduces the first sequence of decimal digits, demonstrating that the largest available dataset of $\pi$ digits can be used to recover $ \pi \approx 3.141593 $ through Monte Carlo simulation.
\end{abstract}

\section{Introduction}

The decimal expansion of $\pi$ has fascinated mathematicians for centuries, not only because of its fundamental role in mathematics and the physical sciences, but also because of its remarkable statistical properties~\cite{BaileyBorwein2016,BorweinBaileyGirgensohn2004}.
Although $\pi$ is a deterministic irrational and transcendental number, its decimal digits appear to behave as if they were generated by a random process. This apparent randomness has motivated extensive theoretical and computational studies aimed at understanding the statistical nature of its decimal expansion.

A central open question concerns the \emph{normality} of $\pi$. A real number is said to be \emph{normal} in base $b$ if every finite sequence of $k$ digits occurs with the expected limiting frequency $b^{-k}$. In base 10, this implies that each digit appears asymptotically with frequency $1/10$, every pair of digits with frequency $1/100$, and, more generally, every block of length $k$ with frequency $10^{-k}$. Although Borel proved that almost all real numbers are normal in the measure-theoretic sense~\cite{Borel1909,KuipersNiederreiter1974,Bugeaud2012}, establishing the normality of specific mathematical constants remains one of the outstanding problems in number theory. In particular, no proof is currently known for $\pi$, despite extensive numerical evidence supporting this conjecture \cite{BaileyCrandall2001, BaileyBorweinPlouffe1997}.

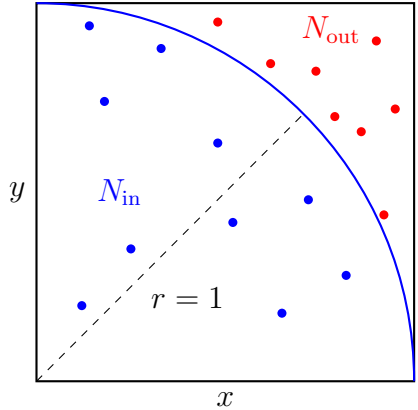
\begin{figure}[t]
\centering
\begin{tikzpicture}[scale=5]

\draw[thick] (0,0) rectangle (1,1);

\draw[thick,blue] (1,0) arc (0:90:1);

\node[below] at (0.5,0) {$x$};
\node[left] at (0,0.5) {$y$};

\foreach \x/\y in {
0.12/0.20,
0.25/0.35,
0.18/0.74,
0.52/0.42,
0.65/0.18,
0.72/0.48,
0.82/0.28,
0.33/0.88,
0.48/0.63,
0.14/0.94
}
\fill[blue] (\x,\y) circle(0.012);

\foreach \x/\y in {
0.95/0.72,
0.86/0.66,
0.92/0.44,
0.74/0.82,
0.62/0.84,
0.48/0.95,
0.90/0.90,
0.79/0.70
}
\fill[red] (\x,\y) circle(0.012);

\node[blue] at (0.22,0.50) {$N_{\rm in}$};
\node[red] at (0.78,0.93) {$N_{\rm out}$};

\draw[dashed] (0,0)--(0.707,0.707);
\node at (0.40,0.22) {$r=1$};

\end{tikzpicture}

\caption{
Illustration of the Monte Carlo estimation of $\pi$.
Random points are uniformly generated inside the unit square.
Points satisfying $x^2+y^2\le1$ (blue) lie inside the quarter unit circle,
whereas the remaining points (red) fall outside.
The ratio of the two areas yields
$\pi \simeq 4N_{\rm in}/N$.
}
\label{fig:mcpi}
\end{figure}

The continuous increase in the number of computed digits of $\pi$ has enabled increasingly stringent statistical tests of its decimal expansion. At the end of 2025, a new computational record reached 314 trillion decimal digits, providing by far the largest dataset ever available for this constant \cite{Yee2025,Backblaze2025}. Statistical analyses performed on smaller dataset showed excellent agreement with the predictions expected for a normal number, although finite computations can never constitute a mathematical proof~\cite{Trueb2016}.

An alternative way to investigate the statistical properties of a sequence is to employ it in stochastic numerical algorithms. Among these, the Monte Carlo estimation of $\pi$ is one of the simplest and most widely known examples~\cite{MetropolisUlam1949,KalosWhitlock2008,Fishman1996}. Consider the unit square $[0,1]\times[0,1]$ containing the quarter of the unit circle centered at the origin, see Fig.~\ref{fig:mcpi}. If points are sampled uniformly within the square, the probability that a point lies inside the quarter circle is equal to the ratio of their areas,

\begin{equation}
P(x^2+y^2\le1)=
\frac{\pi/4}{1}
=\frac{\pi}{4}.
\end{equation}

Consequently, if $N$ random points are generated and $N_{\mathrm{in}}$ fall inside the quarter circle, $\pi$ can be estimated as

\begin{equation}
\hat{\pi}
=
4\,\frac{N_{\mathrm{in}}}{N}.
\end{equation}

The accuracy of this estimator improves as $N^{-1/2}$ from Poissonian fluctuations and depends solely on the availability of uniformly distributed random numbers.

In this work we reverse the conventional approach. Instead of employing a pseudorandom number generator, we use the decimal digits of $\pi$ themselves as the source of random samples for the Monte Carlo algorithm. This approach does not constitute a proof of normality, nor can it establish statistical independence, which cannot hold for a deterministic sequence. Rather, it provides an empirical assessment of the extent to which the known decimal expansion of $\pi$ behaves as a practical source of randomness. By exploiting the complete 314-trillion-digit dataset and optimizing the mapping of decimal digits into Monte Carlo coordinates, we obtain the most accurate self-generated Monte Carlo estimate of $\pi$ currently achievable.

\subsection{The 314-Trillion-Digit Dataset}

The decimal expansion employed in this work originates from the record
computation of $\pi$ completed by StorageReview at the end of 2025 using the
\texttt{y-cruncher} software developed by Alexander J. Yee~\cite{Yee2025}. The computation
reached 314 trillion decimal digits and established a new world record for the
largest contiguous calculation of $\pi$ performed on a single server. The
calculation required approximately 102 days of uninterrupted execution,
demonstrating the maturity of modern high-performance storage architectures as
well as the scalability of state-of-the-art multiple-precision arithmetic
algorithms. 

The computation was performed on a Dell PowerEdge R7725 equipped with two
AMD EPYC 9965 processors, providing a total of 384 CPU cores, 1.5~TB of DDR5
memory, and a storage subsystem consisting of forty 61.44~TB NVMe SSDs for a
total raw capacity of approximately 2.5~PB. Unlike previous record
computations that relied on distributed storage infrastructures, the entire
calculation was executed on a single machine using a carefully optimized local
NVMe array capable of sustaining extremely high I/O bandwidth throughout the
run. The software employed the Chudnovsky series together with the advanced
multiple-precision algorithms implemented in \texttt{y-cruncher}, where the
dominant bottleneck is no longer floating-point performance but rather memory
bandwidth and disk throughput. 

During the computation, intermediate checkpoints occupied more than
1.4~PiB of logical storage, while the total amount of data read and written
exceeded 240~PiB. These figures illustrate that modern record computations of
$\pi$ are essentially large-scale data-management problems, where efficient
I/O scheduling and storage reliability are as important as computational
performance.

Once the calculation had been verified, the final decimal expansion was made
public through the Backblaze B2 cloud storage platform. The published dataset
contains over 130~TB of decimal digits organized into 628 files of
approximately 206~GB each, making it the largest publicly accessible
collection of $\pi$ digits available to date. The much larger temporary
checkpoint files generated during the computation were discarded after
verification, while the final digit sequence was preserved for public access.

Backblaze stores the dataset using its Vault architecture, in which each file
is divided into multiple data shards protected by Reed--Solomon erasure coding
and distributed across independent storage pods. This design provides
high durability while allowing efficient retrieval of very large files,
making the complete decimal expansion readily accessible for mathematical,
statistical, and computational investigations such as the one presented in
this work. 

\section{Construction of the Monte Carlo Sample}

\begin{figure}[t]
\centering
\begin{tikzpicture}[scale=5]

\draw[step=1/30,black!30,very thin] (0,0) grid (1,1);

\draw[thick] (0,0) rectangle (1,1);

\draw[very thick,blue] (1,0) arc (0:90:1);

\foreach \x/\y in {
0.12/0.20,
0.25/0.35,
0.18/0.74,
0.52/0.42,
0.65/0.18,
0.72/0.48,
0.82/0.28,
0.33/0.88,
0.48/0.63,
0.14/0.94,
0.95/0.72,
0.86/0.66,
0.92/0.44,
0.74/0.82,
0.62/0.84,
0.48/0.95,
0.90/0.90,
0.79/0.70
}
\fill ( \x,\y ) circle (0.010);

\draw[dashed] (0,0)--(0.7071,0.7071);
\node at (0.39,0.22) {$r=1$};

\node[below] at (0.5,0) {$x$};
\node[left] at (0,0.5) {$y$};

\end{tikzpicture}

\caption{
Schematic representation of the grid optimization procedure. Consecutive
blocks of decimal digits of $\pi$ are mapped into coordinates $(x,y)$ inside the
unit square. The grid size $N_x\times N_y$ is selected to minimize the probability
of multiple points falling into the same cell (pile-up). For the 314-trillion-digit
dataset, the optimal choice provides a sufficiently large phase space such that
the generated Monte Carlo samples remain statistically independent at the scale
relevant for the estimation of $\pi$.
}
\label{fig:toyMC}
\end{figure}
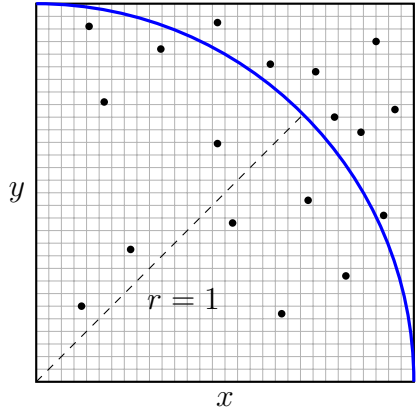

The decimal expansion of $\pi$ is converted into a sequence of two-dimensional
points by partitioning the digit stream into consecutive blocks.
For each point, the first $n_x$ digits are interpreted as the decimal
fraction defining the $x$ coordinate,

\begin{equation}
x=0.d_1d_2\cdots d_{n_x},
\end{equation}
while the following $n_y$ digits define

\begin{equation}
y=0.d_{n_x+1}\cdots d_{n_x+n_y}.
\end{equation}
Successive points are obtained by repeating the same procedure over the
remaining digits. A suitable choice of $n_x$ and $n_y$ is required to avoid an excessive
occupation of the same cells in the $(x,y)$ plane. Let the unit square be
partitioned into a grid of $N_x\times N_y$ cells. Requiring, on average,
one Monte Carlo point per cell gives

\begin{equation}
N_xN_y\simeq 314\times10^{12},
\end{equation}
where $314\times10^{12}$ is the total number of available decimal digits.
Assuming a square grid as reported in Fig.~\ref{fig:toyMC},

\begin{equation}
N_x=N_y=\sqrt{314\times10^{12}}
\simeq1.772\times10^{7}.
\end{equation}
Since each coordinate is encoded by a decimal integer, the required number
of digits is

\begin{equation}
n_x=n_y=\left\lceil\log_{10}N_x\right\rceil=7.
\end{equation}
Therefore each Monte Carlo point consumes $n_x+n_y=14$ decimal digits, yielding a total sample size

\begin{equation}
N=\frac{314\times10^{12}}{14}
=22.4\times10^{12}.
\end{equation}
The statistical uncertainty of the Monte Carlo estimator scales as

\begin{equation}
\sigma_{\pi}\sim\frac{1}{\sqrt{N}}
\sim 10^{-7},
\end{equation}
corresponding to an expected accuracy of approximately seven decimal
digits for the estimate of $\pi$.

\section{Parallel Processing of the Dataset}

The complete decimal expansion of $\pi$ published through the Backblaze cloud
consists of 628 independent files, each containing approximately 206~GB of
decimal digits. Such an organization naturally lends itself to a parallel
processing strategy, since each file can be analyzed independently without
requiring communication with the remaining portions of the dataset.

The multi-threaded simulation code was developed in R with C++ extensions for the binary decoding of YCD files. It was run on an INFN Cloud server (physically hosted at CNAF) with 32 available cores. The data were accessed remotely via rclone, without creating a local copy, resulting in a total execution time of approximately 10 days.
Each worker process reads one chunk of the digit sequence and partitions it
into consecutive blocks of sixteen decimal digits. The first eight digits are
interpreted as the fractional part of the $x$ coordinate,

\begin{equation}
x=0.d_1d_2\cdots d_8,
\end{equation}
while the following eight digits define

\begin{equation}
y=0.d_9d_{10}\cdots d_{16}.
\end{equation}

Each pair therefore represents one point uniformly distributed over the unit
square, assuming that the decimal expansion behaves as a pseudorandom
sequence. The point is classified according to the condition

\begin{equation}
x^2+y^2\le1,
\end{equation}
which determines whether it lies inside the quarter unit circle.
Rather than storing every generated point, each worker maintains only two
integer counters: the total number of processed points, $N_i$, and the number
of points inside the quarter circle, $N_{i,\mathrm{in}}$. At the end of the
analysis of a chunk, the worker returns these two values (and additional data data) to the main thread. Since only
two integers are produced for each input file, the reduction stage requires a
negligible amount of memory and communication.

Each chunk corresponding to one of the 628 samples provides an independent estimate of $\pi$, calculated as:
\begin{equation}
\hat{\pi}_i = 4\frac{N_{\mathrm{in,i}}}{N_i}.
\end{equation}
We performed a one-sample Student's $t$-test to verify that the sample mean is compatible with the theoretical value of $\pi$, obtaining a $p$-value of 0.17. This result indicates that there is no statistically significant evidence to reject the null hypothesis that the mean of the $\hat{\pi}_i$ estimates is equal to the theoretical value of $\pi$.

Under this hypothesis, the maximum likelihood estimate of $\hat{\pi}$ can be obtained by averaging the $\hat{\pi}_i$, since all $N_i$ are nearly identical (with variations at the ppm level). The standard error of $\hat{\pi}$, denoted by $\sigma_{\hat{\pi}}$, can be estimated using the central limit theorem as the standard deviation of the $\hat{\pi}_i$ divided by the square root of the number of chunks
\begin{equation}
\boxed{\hat{\pi} = 3.14159313 \pm 0.00000035}
\end{equation}
reproducing $\pi$ with the expected 7-digit precision, within the estimated uncertainty.

A second approach to estimate both $\hat{\pi}$ and its uncertainty is to exploit the binomial nature of the process using the full statistics, with a success probability $p = \nicefrac{\pi}{4}$, as discussed in Sec. 1. We have
\begin{equation}
\hat{\pi}_B = 4~\frac{\Sigma_i N_{\mathrm{in}, i}}{\Sigma_i N_i} = 4~\frac{17,615,350,257,158}{22,428,557,114,832} = \hat{\pi}
\end{equation}
where the last equality holds only if all $N_i$ are identical. In our data set, $\hat{\pi}_B$ and $\hat{\pi}$ differ by $10^{-13}$, which is entirely negligible. The standard error is
\begin{equation}
\sigma_{\hat{\pi}_B} = 4\sqrt{\frac{p (1 - p)}{\Sigma_i N_i}} = 0.00000035 = \sigma_{\hat{\pi}}.
\end{equation}
The agreement between $\sigma_{\hat{\pi}_B}$ and $\sigma_{\hat{\pi}}$ is expected because each chunk contains a very large number of points, making the sample mean and the full binomial estimator statistically equivalent.

The excellent agreement  ($\Delta=\pi-\hat{\pi}\approx1.4~\sigma_{\hat{\pi}}$) with the theoretical expectation, demonstrating that
the decimal expansion of $\pi$ can successfully act as its own source of pseudorandom numbers in a large-scale Monte Carlo simulation.

The statistical significance of the obtained estimate can be evaluated by
comparing the measured value $\hat{\pi}$ with the exact value of $\pi$ within
the expected fluctuations of the Monte Carlo method. 
Since we have already verified that the standard errors obtained from the central limit theorem and from the sample population are consistent, the $\chi^2$ can be expressed in terms of known quantities:
\begin{equation}
\chi^2=\left(\frac{\Delta}{\sigma_{\hat{\pi}}}\right)^2 = 1.9\, ,
\end{equation}
which follows the known distribution with one degree of freedom under the
hypothesis that the digits of $\pi$ behave as a statistically random source
for the Monte Carlo sampling.
The resulting $p$-value, 

\begin{equation}
P(\chi^2,1)={0.17},
\end{equation}
quantifies the compatibility of the observed deviation with the statistical
fluctuations expected from a Monte Carlo process. A value consistent with the
expected $\chi^2$ distribution indicates that the decimal expansion of $\pi$,
despite being deterministic, provides a statistically valid source for the
generation of Monte Carlo samples at the scale investigated in this work.

Additionally, for each chunk, we filled a histogram with $10^7$ bins, each corresponding to a possible outcome of the 7-digit random numbers. We tested whether the occurrences of the numbers were unbiased. For each histogram, we evaluated the goodness of fit using the Pearson $\chi^2$ test under the hypothesis that the content of each bin was a realization of a Poisson process with mean value $\mu_i = \nicefrac{2 N_i}{10^7}$. The resulting 628 $p$-values are uniformly distributed between 0 and 1, allowing us to exclude biases in the decoding of the binary data or bugs in the analysis code.

\section{Philosophical Perspective}

The present analysis raises a deeper question concerning the relationship
between determinism and randomness in mathematical objects. If the complete
infinite decimal expansion of $\pi$ were available, could $\pi$ itself become a
source capable of generating $\pi$ indefinitely? In other words, could a
deterministic sequence contain enough internal statistical structure to
reconstruct the very object from which it originates? From a strict mathematical point of view, the answer depends on the meaning
assigned to randomness. The digits of $\pi$ are not random in the algorithmic
sense: they are uniquely determined by the definition of the constant.
Consequently, they cannot generate new information beyond what is already
contained in the infinite sequence itself. However, if $\pi$ is normal, its
digits would contain every finite pattern with the correct asymptotic
frequency, making the sequence statistically indistinguishable from a random
source for many practical purposes.

In this hypothetical limit of infinite knowledge, a self-generating process
based on the digits of $\pi$ would not represent a true generation of new
information, but rather a form of self-extraction. The constant would act as
both the encoded object and the apparent source of randomness used to recover
it. This resembles the distinction between randomness and complexity in
algorithmic information theory: a sequence may be completely deterministic
while exhibiting maximal statistical disorder. Therefore, an infinite expansion of $\pi$ could, in principle, generate
arbitrarily accurate Monte Carlo reconstructions of $\pi$ itself, provided that
the sequence possesses the statistical properties expected from a normal
number. The process would not violate determinism, since the output is already
implicitly contained in the input. It would instead reveal a remarkable
property of mathematical information: a sufficiently complex deterministic
object can behave as its own random environment.

\section*{Conclusion}

In this work we have investigated the statistical properties of the largest
currently available decimal expansion of $\pi$ by adopting a complementary
approach to conventional frequency-based tests of normality. Rather than
analyzing the occurrence of digit patterns directly, we employed the decimal
digits themselves as the source of pseudorandom numbers in a Monte Carlo
estimation of $\pi$. A deterministic mapping of consecutive digit blocks into two-dimensional
coordinates was introduced, allowing the construction of nearly
$2\times10^{13}$ independent Monte Carlo samples from the complete
314-trillion-digit dataset. The resulting estimate reproduces the first seven
decimal digits of $\pi$, in agreement with the statistical accuracy expected
from the finite sample size, whose uncertainty scales as $N^{-1/2}$.

Although the present analysis does not constitute a proof of normality, nor
can it establish statistical independence of the digits, it provides a
stringent empirical validation of their random-like behavior in a practical
stochastic computation. In particular, the fact that the decimal expansion of
$\pi$ can successfully serve as its own source of pseudorandom numbers
demonstrates that no detectable statistical bias emerges at the scale of the
largest dataset currently available. This work provides further support for the use of optimized spigot algorithms in the development of new random number generators for Monte Carlo simulations in physics, exploiting arbitrary irrational numbers as seeds.

As computations of $\pi$ continue to extend to even larger numbers of digits,
the same methodology can be naturally refined to produce progressively more
accurate self-generated Monte Carlo estimates. Beyond its intrinsic
mathematical interest, this approach offers an intuitive and quantitative
connection between the concepts of normality, pseudorandomness, and stochastic
simulation, providing an alternative framework for assessing the statistical
properties of one of mathematics' most celebrated constants.

\section*{Acknowledgments}

The authors sincerely thank Brian Beeler (Backblaze) for his support and
assistance in accessing the 314-trillion-digit $\pi$ dataset. The public
availability of this unique computational resource made the present work
possible, and we greatly appreciate his commitment to facilitating its use by
the scientific community.

\end{document}